\documentclass[11pt,reqno]{amsart}
\usepackage{amsbsy,amsfonts,amsmath,amssymb,amscd,amsthm}
\usepackage{mathrsfs,float}
\usepackage[mathcal]{euscript}
\usepackage{enumitem,graphicx,epstopdf}
\usepackage[abs]{overpic}
\usepackage[a4paper,text={6.1in,8in},centering]{geometry}
\usepackage{pxfonts,epstopdf}
\usepackage[colorlinks=true,linkcolor=blue,citecolor=green]{hyperref}
\usepackage{srcltx}
\usepackage[margin=20pt,font=small,labelfont=bf,textfont=it]{caption}
\usepackage{graphicx}
\usepackage{subcaption}
\usepackage{inputenc}
\usepackage{xfrac}

\usepackage{xcolor}
\usepackage[normalem]{ulem}
\usepackage{tikz}
\usepackage{pgfplots}
\pgfplotsset{compat=1.17}

\hypersetup{linktocpage}

\small\normalsize
\theoremstyle{plain}
\newtheorem{mainthm}{Theorem}

\newtheorem{teo}{Theorem}[section]
\newtheorem{claim}[teo]{Claim}
\newtheorem{lem}[teo]{Lemma}
\newtheorem{prop}[teo]{Proposition}

\newtheorem{defi}[teo]{Definition}

\theoremstyle{definition}

\theoremstyle{remark}
\newtheorem{obs}[teo]{Remark}

\makeatletter
\newcommand{\mylabel}[2]{#2\def\@currentlabel{#2}\label{#1}}
\makeatother

\DeclareMathOperator{\IFS}{\mathrm{IFS}}

\newcommand{\eqdef}{\stackrel{\scriptscriptstyle\rm def}{=}}

\setlist[enumerate,1]{label=(\arabic*)}
\setlist[enumerate,2]{label=(\alph*)}

\begin{document}
~\vspace{-1cm}
\title[Lyapunov exponents and IFS]{Zero Lyapunov exponents in transitive skew products of Iterated Function systems}
\begin{abstract}
We study the class of transitive skew products associated with iterated function systems of circle diffeomorphisms. We approximate any of those skew products by maps in this class with a robustly zero Lyapunov exponent. In particular, we prove the existence of non-hyperbolic ergodic measures for an open and dense subset of transitive skew products. Moreover, these measures have full support and are the weak$^*$ limit of periodic measures. 
\end{abstract}

\author[Barrientos]{Pablo G.~Barrientos}
\address{\centerline{Instituto de Matem\'atica e Estat\'istica, UFF}
    \centerline{Rua M\'ario Santos Braga s/n - Campus Valonguinhos, Niter\'oi,  Brazil}}
\email{pgbarrientos@id.uff.br}
\author[Cisneros]{Joel Angel Cisneros}
\address{\centerline{Instituto de Matem\'atica e Estat\'istica, UFF}
    \centerline{Rua M\'ario Santos Braga s/n - Campus Valonguinhos, Niter\'oi,  Brazil}}
\email{joelangel@id.uff.br}

\maketitle  
\thispagestyle{empty}


\section{Introduction}
The theory of uniformly hyperbolic dynamical systems was initiated in the 1960s by Smale~\cite{SMA67} and Anosov~\cite{anosov1967geodesic}, providing a detailed description of a large class of systems, often with very complex evolution. Since the 1970s, with the works of Abraham and Smale~\cite{AS70}, and Newhouse~\cite{NEW70}, it has been known that uniformly hyperbolic systems do not form a dense set of dynamical systems. In response weaker concepts emerged, such as the non-uniformly hyperbolic systems introduced by Pesin~\cite{PES77} in 1977. These systems are characterized by having hyperbolic ergodic measures, that is, with all Lyapunov exponents non-zero. It took a few decades to demonstrate that non-uniformly hyperbolic systems are also not dense. Namely, Kleptsyn and Nalsky~\cite{KN07} constructed in 2007 a dynamical system on the torus $\mathbb{T}^3$ with a robustly zero Lyapunov exponent. Later, in 2016, Bochi, Bonatti, and Díaz~\cite{BBD16} presented new open sets of diffeomorphisms in higher-dimensional manifolds with non-hyperbolic ergodic measures, that is with zero exponents. The result shown by Kleptsyn and Nalsky is based on a toy model built a few years earlier by Gorodetski, Ilyashenko, Kleptsyn, Nalsky~\cite{GI}. In this paper, we will show that the existence of non-hyperbolic ergodic measures is abundant in the class of transitive skew products associated with iterated function systems where the toy models in~\cite{GI} were constructed. 

An iterated function system (IFS) on $\mathbb{S}^1$ is a family of functions $\mathscr{F}=\{f_1,\dots,f_k\}$ from the circle~$\mathbb{S}^1$ to itself that can be applied (composed) successively in any order. Therefore, such compositions are elements of the semigroup $\langle \mathscr{F} \rangle^+$ associated with these transformations. Moreover, the sequence of compositions can be seen as the fiber of the iterations of the locally constant (or one-step) skew product 
\begin{equation} \label{eq:skew product}
  F:\Sigma_k\times \mathbb{S}^1 \to \Sigma_k \times \mathbb{S}^1, \qquad
  F^n(\omega,x)=(\sigma^n(\omega),f^n_\omega(x))
\end{equation}
where at the base we have the shift $\sigma:\Sigma_k\to \Sigma_k$ and
$$
   f^n_\omega\eqdef f_{\omega_{n-1}}\circ \dots \circ
   f_{\omega_0} \quad \text{for} \ \
   \omega=(\omega_i)_{i\geq 0} \in \Sigma_k\eqdef
   \{1,\dots,k\}^{\mathbb{N}} \quad \text{and} \quad n\geq 1.
$$ 
To emphasize the role of the fiber maps and the product structure, we will write $F=\sigma\ltimes \mathscr{F}$.



Given an ergodic $F$-invariant measure $\mu$, Birkhoff's Ergodic Theorem guarantees the existence of the limit
$$
   \lambda(\omega,x)=\lim_{n\to \infty} \frac{1}{n}\log |(f^n_\omega)'(x)| \quad \text{for $\mu$-a.e.~$(\omega, x)\in \Sigma_k\times
   \mathbb{S}^1$}.
$$
In fact
\begin{equation} \label{eq:lyp}
   \lambda(\mu)\eqdef \int \log |(f_{\omega})'(x)| \, d\mu
   =\lambda(\omega,x) \quad  \text{for $\mu$-a.e   $(\omega,x)\in \Sigma_k\times
   \mathbb{S}^1$}.
\end{equation}
The value $\lambda(\mu)$ is called the \emph{Lyapunov exponent} along the fiber of the skew product $F$  
with respect to the measure $\mu$.

For $r\geq 1$ and $k\geq 2$, we  denote by $\IFS_k^r(\mathbb{S}^1)$ the set of families $\mathscr{F}=\{f_1,\dots,f_k\}$ with $f_i$ belongs to the space $\mathrm{Diff}^r_+(\mathbb{S}^1)$ of orientation-preserving diffeomorphisms of the circle equipped with the $C^r$-topology. Moreover, we say that $\mathscr{F}=\{f_1,\dots,f_k\}$ and  $\mathscr{G}=\{g_1,\dots,g_k\}$ are $C^r$-close if, by reordering the elements families if necessary, $f_i$ and $g_i$ are $C^r$-close for every $i=1,\dots,k$.  We define $\mathrm{S}_{k}^{r}(\mathbb{S}^{1})$ as the set of skew products $F=\sigma\ltimes \mathscr{F}$, where $\mathscr{F}\in \IFS_k^r(\mathbb{S}^1)$. We endow $\mathrm{S}_{k}^{r}(\mathbb{S}^{1})$ with a topology, stating that $F=\sigma\ltimes \mathscr{F}$ is close to $G=\sigma\ltimes \mathscr{G}$ if $\mathscr{F}$ and $\mathscr{G}$ are $C^r$-close.

Gorodetski et al.~\cite{GI} constructed an open set $\mathcal{U}$ in $\mathrm{S}_2^1(\mathbb{S}^1)$
such that  each  ${F}\in \mathcal{U}$
has an  ergodic invariant measure $\mu$ with $\lambda(\mu)=0$.  Such measures are approximated by periodic measures with a specific tail structure and called \emph{$GIKN$-measures} or \emph{royal measures}, c.f.~\cite{KL17}.
Originally in~\cite{GI}, it is shown that this open set can be constructed as a $C^r$-neighborhood of a skew product  $G=\sigma\ltimes \mathscr{G}$ with $\mathscr{G}=\{g_1,g_2\}$ where $g_1$ is  an irrational rotation and $g_2$ has a pair of hyperbolic periodic points with additional conditions on the eigenvalues. Later, Bochi et al.~\cite{BBD13} generalized this result for IFSs on compact manifolds of higher dimension. In fact, they constructed open sets of IFSs in the $C^r$ topology for $r\geq2$ where the associated skew products have an ergodic invariant measure, limit in the weak$^*$ topology of periodic measures, with full support and where all Lyapunov exponents along the fibers are zero. 
It is worth noting that the measures obtained by the method of approximation by periodic orbits typically have zero entropy, as proved in \cite{KL17}. Further developments regarding the construction of non-hyperbolic measures include the $GIKN$ method based on so-called \emph{flip-flop families} introduced in \cite{BBD16, BDB18}, which provides a set of positive entropy supporting only non-hyperbolic measures. More recently, ``mixed'' methods combining the $GIKN$ and flip-flop approaches have established the robust existence of non-hyperbolic ergodic measures with both positive entropy and full support \cite{BDK21, Lac23}. For a detailed survey of these results, we refer the reader to \cite{Dia18}.





A skew product $F\in \mathrm{S}_{k}^{r}(\mathbb{S}^{1})$ is said to be \emph{transitive} if there has a dense orbit. 
We denote by $\mathrm{TS}_k^r(\mathbb{S}^1)$ the subset of $\mathrm{S}_{k}^{r}(\mathbb{S}^{1})$ formed by transitive skew products. The open set of skew products constructed in~\cite{GI} is contained in this set. The following main result shows that any transitive skew product can be approximated by an open  (in $\mathrm{S}^r_k(\mathbb{S}^1)$) set of transitive skew products having a non-hyperbolic ergodic measure with full support. This concludes affirmatively the conjecture~\cite[Conj.~2]{BD} in our setting.

\begin{mainthm} \label{thmAA}
For $k\geq 2$ and $1\leq r\leq \infty$, there exists an open and dense  subset $\mathcal{R}$ of $\mathrm{TS}_{k}^{r}(\mathbb{S}^{1})$ such that for every $F\in\mathcal{R}$ there is an ergodic $F$-invariant probability measure $\mu$ with full support whose Lyapunov exponent $\lambda(\mu)$ is equal to zero. The probability measure $\mu$ is a $GIKN$-measure.
\end{mainthm}

The $C^1$-density of maps with non-hyperbolic ergodic invariant measures, as in the above theorem, was stated in~\cite[Cor.~2 and~8.10]{diaz2017nonhyperbolic} on the set of all \emph{locally constant}\footnote{Although the result in~\cite{diaz2017nonhyperbolic} is stated for H\"older skew products, the approximation by locally constant skew products in cylinders proposed is not clear.} skew products with $C^1$-fiber maps that are robustly transitive and have periodic points of different indices.
This result is a corollary of Proposition~8.9 in~\cite{diaz2017nonhyperbolic}, which is not fully proved. The authors observe that the proof of a similar result in~\cite{bonatti2002minimality} for partially hyperbolic diffeomorphisms of compact manifolds can be translated \emph{mutatis mutandis} to obtain such a proposition in the locally constant skew product setting. The results in~\cite{bonatti2002minimality} show the density of the strong unstable (strong stable) foliation, and the authors in~\cite{diaz2017nonhyperbolic} translate this property to the locally constant skew product setting as the forward (backward) minimality of the underlying IFS. However, this translation is incorrect: while the minimality of the strong unstable foliation implies the forward minimality of the IFS, the converse does not hold. As demonstrated by the counterexamples in~\cite[Cor.~6.10]{BC23}, there exist IFSs that are forward minimal, yet their associated strong unstable foliations are not minimal. Consequently, establishing the density of minimal IFSs is insufficient to deduce the density of minimal strong foliations, leaving the proof in~\cite{diaz2017nonhyperbolic} incomplete. Nevertheless, according to~\cite[Thm.~B]{BC23}, it follows that the set of skew products having both strong stable and strong unstable foliations also contains an open and dense subset of $\mathrm{TS}^r_k(\mathbb{S}^1)$. Consequently, the set $\mathcal{R}$ in Theorem~\ref{thmAA} not only satisfies the conclusions of the theorem but can also be chosen to include this additional geometric property.

In \S\ref{s-thmB}, we revisit the proof of a key result from Gorodetski et al.~\cite{GI}, which introduced the method of approximation by periodic orbits. This method serves as a fundamental tool for our main theorem. Later works~\cite{Diaz2009non, BD, BBD13} refined and extended this approach. By incorporating these refinements, we simplify and improve the original proof given in~\cite{GI}. In~\S\ref{s-thmA}, we use this framework and recent results from~\cite{BC23} to prove Theorem~\ref{thmAA}.  

\section{Ergodic non-hyperbolic measure with full support}  \label{s-thmB}

To prove Theorem~\ref{thmAA} we need to the following slight generalization of~\cite[Thm.~2]{GI}. First, we recall some definitions. 

A finite family $\mathscr{F}$ of circle diffeomorphisms is
\begin{enumerate}[label=--]
    \item \emph{minimal} if for every $x\in \mathbb{S}^1$ the orbit $\mathcal{O}(x)=\{f(x): f\in \langle \mathscr{F}\rangle^+ \}$ is a dense set in~$\mathbb{S}^1$;
    \item \emph{expanding} if for every $x\in \mathbb{S}^1$, there exists $f$ in $\langle \mathscr{F}\rangle^+$ such that $|(f^{-1})'(x)|>1$. 
\end{enumerate}
We say that $\mathscr{F}$ is \emph{backward minimal} and \emph{backward expanding} if the family $\mathscr{F}^{-1}=\{f^{-1}: f\in \mathscr{F}\}$ is minimal and expanding respectively.  


\begin{teo} \label{teo-B}
For any $k\geq 2$, consider the skew product $F=\sigma\ltimes \mathscr{F}$ where  $\mathscr{F}=\{f_1,\dots,f_k\}$ is a family of circle $C^1$-diffeomorphisms such that  
\begin{enumerate}[label=(\roman*)]
  \item \label{item1-teo-B}$\mathscr{F}$ is minimal;
  \item \label{item2-teo-B}$\mathscr{F}$ is backward expanding; 
  \item \label{item3-teo-B}there exists $f\in \langle\mathscr{F}\rangle^+$
  with a hyperbolic attracting periodic point.
\end{enumerate}
Then  $F$  has an ergodic invariant measure $\mu$ with 
 full support, which is a $GIKN$- measure and satisfies $\lambda(\mu)=0$.
\end{teo}

The method of approximation by periodic orbits, introduced in~\cite{GI}, is a central tool for constructing the non-hyperbolic invariant measures in the above theorem. However, the property of full support of the limit measure was not directly proven in~\cite{GI}.  This property was obtained later in~\cite{BBD13} for skew products with higher dimensional fiber space by spreading the sequence of periodic points in the ambient space using previous results from~\cite{BD}. In what follows, we replicate these ideas to demonstrate the above~theorem.

\subsection{Preliminaries}

In this section, we collect some results from~\cite{GI, Diaz2009non, BD, BBD13, BFMS16} that will be used in the proof of Theorem~\ref{teo-B}. 


The following proposition, which follows directly from the definition of weak$^*$ topology and~\eqref{eq:lyp}, is fundamental to the Lyapunov exponent approximation method.

\begin{prop}\label{pro3.3-tesis}
Let $\{\mu_n\}_{n \in \mathbb{N}}$ and $\mu$ be invariant ergodic probability measures for the skew product $F = \sigma \ltimes \mathscr{F}$. 
If $\mu_n \to \mu$ as $n \to \infty$ in the weak$^*$ topology, then 
$\lambda(\mu_n) \to \lambda(\mu)$ as $n \to \infty$.
\end{prop}

Given an $n$-periodic orbit 
\[
X = \left\{ G(x_0), G^2(x_0), \dots, G^{n-1}(x_0), G^n(x_0) = x_0 \right\}
\]
of a continuous map $G$ on a compact metric space, the atomic measure uniformly distributed on $X$ is defined by  
\[
\mu_X = \frac{1}{n} \sum_{i=0}^{n-1} \delta_{G^i(x_0)},
\]
where $\delta_y$ denotes the Dirac measure at the point $y$. This measure is $G$-invariant and ergodic.

\begin{defi}\label{defi-boa-aprox}
Let $X$ and $Y$ be periodic orbits of $G$. Given $\gamma > 0$ and $\aleph > 0$, we say that $Y$ is a \emph{$(\gamma, \aleph)$-good approximation} of $X$ if the following conditions hold:
\begin{enumerate}[label=(\roman*)]
    \item \label{item-boa1} There exists a subset $\Gamma \subset Y$ and a projection $\rho : \Gamma \to X$ such that $d(G^j(y), G^j(\rho(y))) < \gamma$ for all $y \in \Gamma$ and all $j = 0, 1, \dots, P-1$, where $P$ is the period of $X$;
    \item \label{item-boa2} $\frac{\#\Gamma}{\#Y} \geq \aleph$;
    \item \label{item-boa3}$\#\rho^{-1}(x)$ is the same for all $x\in X$.
\end{enumerate}
\end{defi}

\begin{defi}
We refer to a sequence $\{X_n\}_{n \in \mathbb{N}}$ of $G$-periodic orbits, where $\#X_n\to \infty$ as $n\to \infty$, as a \emph{$GIKN$-sequence} 
if there exist sequences of positive real numbers $\{\gamma_n\}_{n\in\mathbb{N}}$ and $\{\aleph_n\}_{n\in\mathbb{N}}$ that satisfy:
\begin{enumerate}
\item\label{item1teoBonati} For each $n \in \mathbb{N}$, the orbit $X_{n+1}$ is a $(\gamma_n, \aleph_n)$-good approximation of $X_n$;
    \item\label{item2teoBonati} $\sum_{n=1}^\infty \gamma_n < \infty$;
    \item\label{item3teoBonati} $\prod_{n=1}^\infty \aleph_n \in (0, 1]$.
\end{enumerate}
By an abuse of terminology, we will also refer to the sequence of atomic measures $\{\mu_n\}_{n\in\mathbb{N}}$ uniformly distributed over the orbit $X_n$, as to \emph{$GIKN$-sequence of measures}.
\end{defi}
The next result corresponds to~\cite[Lem.~2.5]{BD}, which refines~\cite[Prop.~2.5]{Diaz2009non}. Compare also with~\cite[Lem.~2 and Section 8]{GI}.

\begin{teo}\label{teoBonati}
Let $\{\mu_n\}_{n\in\mathbb{N}}$ be the $GIKN$-sequence of ergodic measures associated with a sequence of $G$-periodic orbits $\{X_n\}_{n \in \mathbb{N}}$. 
Then $\{\mu_n\}_{n \in \mathbb{N}}$  weak$^*$ converges to an ergodic measure $\mu$  supported on the topological limit of $\{X_n\}_{n \in \mathbb{N}}$, that is,
\[
\mathrm{supp}~\mu = \bigcap_{k=1}^\infty \overline{\bigcup_{l=k}^\infty X_l}.
\]
\end{teo}

\begin{defi}
A measure $\mu$ is \emph{$GIKN$-measure} if it is a weak$^*$ limit of a $GIKN$- sequence.
\end{defi}

Given a finite word $\bar{\omega} = \bar{\omega}_0 \bar{\omega}_1 \dots \bar{\omega}_{P-1}$ in the alphabet $\{1, \dots, k\}$, i.e., a sequence of length $|\bar{\omega}| = P$ with $\bar{\omega}_j \in \{1, \dots, k\}$ for $j = 0, 1, \dots, P-1$, we define the cylinder 
\[
\llbracket \bar{\omega} \rrbracket \eqdef \{ \omega = (\omega_i)_{i \geq 0} \in \Sigma_k : \omega_i = \bar{\omega}_i \text{ for } i = 0, \dots, |\bar{\omega}|-1 \},
\]
the composition
\[
T_{\bar{\omega}} \eqdef f_{\bar{\omega}_{P-1}} \circ f_{\bar{\omega}_{P-2}} \circ \dots \circ f_{\bar{\omega}_0},
\]
and use $\bar{\omega}\,\overset{n}{\ldots}\,\bar{\omega}$ to represent the concatenation of $n$ copies of $\bar{\omega}$. Furthermore, the infinite concatenation of copies of $\bar{\omega}$ is denoted by $\bar{\omega}^\infty$.

\begin{obs}\label{obs-f-gera-orb-perio}
Clearly, if $f\in\langle \mathscr{F} \rangle^+$, then there is an integer $P\geq 1$ and a finite word $\bar{\omega}=\bar{\omega}_0 \bar{\omega}_1 \dots \bar{\omega}_{P-1}$ such that $f=T_{\bar{\omega}}$. Furthermore, if $f$ has a fixed point $x\in\mathbb{S}^1$, then $z=(\bar{\omega}^\infty,x)\in\Sigma_k\times\mathbb{S}^1$ is a $P$-periodic point of the skew product $F = \sigma \ltimes \mathscr{F}$. Hence, $f$ generates a $P$-periodic orbit $X=\{z,F(z),\dots,F^{P-1}(z)\}$ of the skew product $F$, where $P$ is the number of maps in $\mathscr{F}$ generating $f$.
\end{obs}

One advantage of using an atomic periodic measure for the skew product $F = \sigma \ltimes \mathscr{F}$ is that its Lyapunov exponent can be computed explicitly. This leads to the following proposition:

\begin{prop}\label{propo-expoente-orbita-perio}
Consider $f \in \langle \mathscr{F} \rangle^+$ with a fixed point $x \in \mathbb{S}^1$. Let $X$ be the $P$-periodic orbit of the skew product $F = \sigma \ltimes \mathscr{F}$ generated by $f$. Then, 
\[
\lambda(\mu_X) = \frac{1}{P} \log |f'(x)|,
\]
where $\mu_X$ denotes the atomic measure uniformly distributed over $X$.
\end{prop}


The following lemma from~\cite[Lem.~6.3]{BBD13}, provides a key tool for establishing density properties of orbits.

\begin{lem}\label{lema-ir-a-casa-dado-um-paseo}
Assume that $\mathscr{F}$ is minimal. Then, for any $\eta > 0$ and any non-empty open set $J \subset \mathbb{S}^1$, there exist constants $\delta_0 = \delta_0(\eta, J) > 0$ and $K = K(\eta, J) \in \mathbb{N}$ such that for every open interval $I \subset \mathbb{S}^1$ with length at most $\delta_0$, there exists a finite word $\bar{\omega}$ in the alphabet $\{1, 2, \dots, k\}$ with $|\bar{\omega}| = R \leq K$ satisfying:
\begin{enumerate}[label={(\Roman*)}, itemsep=0.2cm]
    \item\label{ec.densidade-da-orbita} For every $(\omega, x) \in \llbracket \bar{\omega} \rrbracket \times I$, the segment of the orbit of $F = \sigma \ltimes \mathscr{F}$,
    \[
    F^{[0,R]}(\omega, x) \eqdef \{ F^j(\omega, x) : 0 \leq j \leq R \},
    \]
    is $\eta$-dense in $\Sigma_k \times \mathbb{S}^1$;
    \item\label{iii.12} $T_{\bar{\omega}}(I) \subset J$.
\end{enumerate}
\end{lem}

The following proposition, which appears in~\cite[Prop.~3]{GI}, provides a complementary result about the control of intervals under compositions of the generators of $\mathscr{F}$.

\begin{prop}\label{propo3-Gorodesky}
If $\mathscr{F}$ is minimal, then for any interval $J$ in $\mathbb{S}^1$, there exist constants $K_0 = K_0(J) \in \mathbb{N}$ and $\delta_0 = \delta_0(J) > 0$ such that any interval $I$ with length smaller than $\delta_0$ can be mapped into $J$ by some composition of at least $K_0$ generators of $\mathscr{F}$.
\end{prop}
\begin{proof}
Assume Lemma~\ref{lema-ir-a-casa-dado-um-paseo} holds. Take any $\eta>0$. Fix $J\subset\mathbb{S}^1$. Use Lemma~\ref{lema-ir-a-casa-dado-um-paseo} to find $\delta_{0}=\delta_{0}(\eta,J)$ and $K=K(\eta,J)$. Cover $\mathbb{S}^1$ with finitely many intervals of length $\delta_{0}/{2}$. Denote the interval forming that cover by $I_{1},\cdots,I_{t}$. Applying Lemma~\ref{lema-ir-a-casa-dado-um-paseo} to each $1\leq \iota\leq t$ we find a word $\omega(\iota)$ with $|\omega(\iota)|\leq K$ such that $T_{\omega(\iota)}(I_{\iota})\subset J$. Let $K'_{0}=\min\{|\omega(\iota)|: 1\leq\iota\leq t\}$. Let $\delta'_{0}$ be the Lebesgue number for our cover. Now consider any interval $I\subset\mathbb{S}^1$ with length less that $\delta'_{0}$. There is $1\leq\iota\leq t$ such that $I\subset I_{\iota}$ and $I$ is mapped into $J$ by a composition of $|\omega(\iota)|\geq K'_{0}$ generators of $\mathscr{F}$. Hence $\delta'_{0}$ and $K'_{0}$ are constants needed to establish Proposition~\ref{propo3-Gorodesky}.
\end{proof}

Finally, the following proposition is an immediate consequence of~\cite[Lem.~3.1]{BFMS16} and provides key properties of backward expanding systems.

\begin{prop}\label{propo4-Gorodesky}
Assume that $\mathscr{F}$ is backward expanding. Then there exist constants $\nu > 1$, functions $h_1, \dots, h_m \in \langle \mathscr{F} \rangle^+$, and open intervals $B_1, \dots, B_m$ of $\mathbb{S}^1$ such that:
\begin{enumerate}[label=(\alph*)]
    \item \label{item-a-consecuencia-propo2.3} $\mathbb{S}^1 = B_1 \cup B_2 \cup \dots \cup B_m$;
    \item \label{item-b-consecuencia-propo2.3} $|h'_i(x)| > \nu$ for all $x \in B_i$ and $i = 1, 2, \dots, m$.
\end{enumerate}
\end{prop}

\subsection{Main Lemma}
This lemma is an improved version of~\cite[Lemma~3]{GI} where we spread the constructed periodic orbits densely on $\Sigma_k\times \mathbb{S}^1$. 

\begin{lem}\label{lema3.1tesis}
Consider $F=\sigma\ltimes \mathscr{F}$ and
assume  that
\begin{enumerate}[label=(\roman*)]
\item\label{item-i-lema3.1} $\mathscr{F}$ is minimal ;
\item\label{item-ii-lema3.1} $\mathscr{F}$ is backward expanding;
\item\label{item-iii-lema3.1} there exists $f\in \langle\mathscr{F}\rangle^+$ with a hyperbolic attracting periodic point.
\end{enumerate}
Then, there exist $0<c<1$ and $d>0$ depending only on $\mathscr{F}$, such that if $X$ is the $P$-periodic orbit of $F$ generated by $f$, then for every $\gamma>0$ and $\eta>0$, there exists a periodic orbit $X'$ of $F$ with period $P'>P$ with the following properties:
\begin{enumerate}
  \item\label{item1-lema3.1} $ c \cdot \lambda(\mu_{X})< \lambda(\mu_{X'}) <0$;
  \item\label{item2-lema3.1} $X'$ is $(\gamma,\aleph)$-good approximation\footnote{Note that a priori, $\aleph$ may not be positive (requirement of the definition of good approximation). However, the lemma provides a new periodic orbit with smaller exponent. Thus, recursively applying the lemma, we always find an element $f$ satisfying (iii) whose associated periodic orbit has an exponent small enough with~$\aleph>0$.} of $X$ where $\aleph\eqdef 1-d\cdot |\lambda(\mu_{X})|$;
  \item\label{item3-lema3.1} $X'$ is $\eta$-dense in $\Sigma_k\times\mathbb{S}^1$.
\end{enumerate}
\end{lem}
\begin{proof}
First, let us prove the existence of the periodic orbit $X'$ of $F$ with a period greater than~$P$. Since  $f \in \langle\mathscr{F} \rangle^{+}$, there exists a finite word $\hat{\omega}$ such that $f = T_{\hat{\omega}}$ and  $P=|\hat{\omega}|$. Let $x$ be the attracting periodic point for $f$, meaning $f(x) = x $ and $|f'(x)| = \alpha $ where $ 0 < \alpha < 1 $. Let us take constants $\alpha_{-} $ and 
$ \alpha_{+} $ such that
\begin{equation}\label{iii.9}
  0 < \alpha_{-} < \alpha < \alpha_{+} < 1.
\end{equation}
Fix $\gamma>0$. By the continuity of $f'$, there is an interval $J $ containing $x$ such that 
\begin{equation}\label{iii.10}
  0 < \alpha_{-} \leq |f'(y)| \leq \alpha_{+} < 1 \quad \text{for all $y\in J$}, \quad\text{and}\quad L^{2P}|J| < \gamma,
\end{equation}
where we denote by $|I|$ the length of an interval $I\subset \mathbb{S}^{1}$ and $L\eqdef\max\{\|g'\|_\infty: g\in \mathscr{F}\}$.

From \ref{item-i-lema3.1} and Lemma~\ref{lema-ir-a-casa-dado-um-paseo}, for $J$ and $\eta > 0$, 
there exist  $\delta_{0} = \delta_{0}(J,\eta) > 0$ and  $K = K(J,\eta) \in \mathbb{N}$ such that~\ref{ec.densidade-da-orbita} and~\ref{iii.12} hold.  
On the other hand, from~\ref{item-ii-lema3.1} and Proposition~\ref{propo4-Gorodesky}, 
 there exist a constant $\nu>1$, diffeomorphisms $h_{1},\dots,h_{m}\in\langle \mathscr{F}\rangle^+$, and open intervals $B_{1},\dots,B_{m}$ of $\mathbb{S}^{1}$ such that~\ref{item-a-consecuencia-propo2.3} and~\ref{item-b-consecuencia-propo2.3} hold. 
Let $\lambda > 0$ be the Lebesgue number of the cover in~\ref{item-a-consecuencia-propo2.3}. Recall that this means any interval 
of length less than $\lambda$ is contained in some~$B_i$. Let us set
\begin{equation}\label{iii.15}
L_{1}=\max\{\|h'_{i}\|_{\infty}: i=1,\dots,m\}>1 \quad \text{and} \quad
\delta=\min\{\delta_{0},|J|L^{-K},\lambda/L_1\}>0.
\end{equation}
For each $n \in \mathbb{N}$ sufficiently large, we take $r = r(n) \in \mathbb{N}$ to be the unique integer such that
\begin{equation} \label{eq:Pablo}
   0 < \frac{\log \frac{\alpha_+^{-n} \delta}{L_1|J|}}{\log L_1} \leq r < \frac{\log \frac{\alpha_+^{-n} \delta}{L_1|J|}}{\log
   L_1} + 1 = \frac{\log\frac{\alpha_+^{-n}\delta}{|J|}}{\log L_1}.
\end{equation}
Note that it holds
\begin{equation}\label{iii.19}
  \frac{\delta}{L_{1}}\leq\alpha_{+}^nL_{1}^{r}|J|<\delta.
\end{equation}
Since  $x \in J$ is an attracting fixed point of $f$, using the Intermediate Value Theorem, the condition $L_{1}>1$, and~\eqref{iii.10},~\eqref{iii.19} and~\eqref{iii.15},  we have
$ |f^n(J)|\leq \alpha_+^n |J| \leq \alpha^n_+ L_1^{r} |J| <
   \delta  < \lambda$.
Thus, there exists $i_1\in\{1,\dots,m\}$ such that $f^n(J)\subset B_{i_1} $. Similarly, since
$$|h_{i_1}f^n(J)|\leq \alpha_+^n
L_1|J|\leq \alpha_+^n L_1^{r}|J|<\delta < \lambda,
$$
we obtain $i_2\in \{1,\dots,m\}$ such that
$h_{i_1}f^n(J)\subset B_{i_2}$. Proceeding recursively,
we obtain $i_1,\dots,i_r\in \{1,\dots,m\}$ such that
\begin{equation} \label{iii.13}
  f^n(J)\subset B_{i_1} \ \ \text{and} \ \  h_{i_j}\dots h_{i_1}f^n(J)\subset B_{i_{j+1}}
  \
    \text{for every  $j=1,\dots,r-1$.}
\end{equation}
Denote by $\bar{\kappa}_0=\hat{\omega}\,\overset{n}{\ldots}\,\hat{\omega}$ and consider the finite word $\bar{\kappa}_1$ such that
$ T_{\bar{\kappa}_1}=h_{i_r}\circ \dots \circ h_{i_1}$.
Then, using again the Intermediate Value Theorem, observing that $T_{\bar\kappa_0}=f^n$ and due to~\eqref{iii.13},~\eqref{iii.10}, and~\ref{item-b-consecuencia-propo2.3}, we have
\begin{equation}\label{iii.20}
 \alpha^n_- \nu^r|J| \leq |T_{\bar\kappa_1}T_{\bar\kappa_0}(J)|\leq L_{1}^{r}\alpha_{+}^{n}|J|.
\end{equation}
Let $I=T_{\bar\kappa_1}T_{\bar\kappa_0}(J)$. Then, from equations~\eqref{iii.19}, \eqref{iii.20}, it follows that $|I|<\delta \leq\delta_0$. Thus, by the definition of $\delta_{0}=\delta_{0}(J,\eta)$ and~\ref{iii.12}, there exists a finite word $\bar\kappa_2$, with $|\bar\kappa_2|\leq K$, such that $T_{\bar\kappa_2}(I)\subset J$. Therefore,  $T_{\bar\kappa_2}T_{\bar\kappa_1}T_{\bar\kappa_0}(J)\subset J$  by the definition of $I$.
Taking the finite word $\bar{\omega}=\bar\kappa_0\bar\kappa_1\bar\kappa_2$, we have $T_{\bar{\omega}}=T_{\bar\kappa_2}\circ T_{\bar{\kappa}_1}\circ T_{\bar\kappa_0}$, and thus,
\begin{equation}\label{iii.22}
  T_{\bar{\omega}}(J)\subset J  \quad \text{with} \quad |\bar{\omega}|=nP+|\bar\kappa_1|+|\bar\kappa_2|>P.
\end{equation}
From $|\bar\kappa_2|\leq K$, \eqref{iii.15}, and \eqref{iii.19}, we have
\begin{equation}\label{iii.27}
 |T'_{\bar{\omega}}(y)|\leq L^KL_1^r\alpha_+^n<
 L^{K} \frac{\delta}{|J|}\leq L^{K}\frac{|J|L^{-K}}{|J|}=1
 \quad \text{for all} \quad y\in J.
\end{equation}
Hence, \eqref{iii.22} and \eqref{iii.27} imply that $T_{\bar{\omega}}$ is a contraction on $J$. Therefore, by the Banach Fixed-Point Theorem, $T_{\bar{\omega}}$ has a unique fixed point $x'\in J$. Finally, by Remark~\ref{obs-f-gera-orb-perio}, there exists a periodic orbit $X'$ generated by $T_{\bar{\omega}}$ with period $P'=|\bar{\omega}|>P$.

It should be noted that the previous argument holds for any sufficiently large integer $n$ and constants $\alpha_+$ and $\alpha_-$ satisfying~\eqref{iii.9}. Furthermore, since $\bar\kappa_2$ is the finite word given by Lemma~\ref{lema-ir-a-casa-dado-um-paseo} and the fact that $\sigma^{nP+\bar{\kappa}_1}(\bar{\omega}^\infty)\in \llbracket\bar\kappa_2\rrbracket$, it follows from~\ref{ec.densidade-da-orbita} that $$F^{[0,|\bar\kappa_2|]}\left(\sigma^{nP+\bar{\kappa}_1}(\bar{\omega}^\infty), T_{\bar\kappa_1}T_{\bar\kappa_0}(x')\right),$$ is $\eta$-dense in $\Sigma_k\times\mathbb{S}^1$. Since $F^{[0,|\bar\kappa_2|]}\left(\sigma^{nP+\bar{\kappa}_1}(\bar{\omega}^\infty), T_{\bar\kappa_1}T_{\bar\kappa_0}(x')\right)\subset X'$, it follows that the orbit $X'$ is $\eta$-dense in $\Sigma_k\times\mathbb{S}^1$, thus proving Item~\ref{item3-lema3.1}.

Now let us calculate the Lyapunov exponent along the fiber with respect to the measure $\mu_{X'}$. For this, according to~\eqref{iii.20}, we have
$
  |T'_{\bar{\omega}}(x')|\geq
  |T'_{\bar\kappa_2}(T_{\bar{\kappa}_1}T_{\bar\kappa_0}(x'))|\,
  \alpha_-^n\nu^r \geq e^{C_1} \alpha_-^n \nu^r
$,
where
\begin{equation}\label{eq:C1}
C_1=\min \{ \log |DT_{\bar\kappa}(y)|: \, y\in \mathbb{S}^1 \ \text{and} \ |\bar\kappa|<K \}.
\end{equation}
Thus, we have
\begin{equation}\label{iii.33}
  \log|T'_{\bar{\omega}}(x')|> C_{1}+r\log\nu+n\log\alpha_{-}.
\end{equation}
Substituting~\eqref{eq:Pablo} into~\eqref{iii.33}, we get
\begin{equation}\label{iii.35}
\begin{aligned}
  \log|T'_{\bar{\omega}}(x')| > C_{1}+ \frac{\log \nu}{\log L_1}
 \log \frac{\alpha_+^{-n}\delta}{L_1 |J|} +n\log\alpha_{-}
  = C_{1}+ C_2 + n\, ( \log \alpha_- - \frac{\log \nu}{\log L_1}\log
  \alpha_+)
\end{aligned}
\end{equation}
where
\begin{equation}\label{eq:C2}
C_2=\frac{\log \nu}{\log L_1} \log \frac{\delta}{L_1 |J|}.
\end{equation}
Taking into account that $0<\alpha<1$ and $1<\nu<L_1$, 
$$
(1-\frac{2\log\nu}{3\log L_{1}})\log \alpha <
(1-\frac{\log\nu}{\log L_{1}})\log \alpha <0.
$$
Thus, since
$$
  \lim_{(\alpha_-,\alpha_+)\to (\alpha,\alpha)} \log\alpha_{-}-\frac{\log\nu}{\log
  L_{1}}\log\alpha_{+} = (1-\frac{\log\nu}{\log L_{1}})\log
  \alpha,
$$
there exist $\alpha_{-}$ and $\alpha_{+}$ satisfying~\eqref{iii.9} and~\eqref{iii.10}, such that
\begin{equation}\label{iii.36}
  (1-\frac{2\log\nu}{3\log L_{1}})\log\alpha \leq \log\alpha_{-}-\frac{\log\nu}{\log
  L_{1}}\log\alpha_{+}.
\end{equation}
Writing $C_{3}=C_{1}+C_{2}$ and substituting~\eqref{iii.36} into~\eqref{iii.35}, it follows that
\begin{equation}\label{iii.37}
  \log|T'_{\bar{\omega}}(x')|>C_{3}+n(1-\frac{2\log\nu}{3\log L_{1}})\log\alpha.
\end{equation}
From Proposition~\ref{propo-expoente-orbita-perio}, the Lyapunov exponent along the fiber with respect to the measure $\mu_{X'}$ can be estimated as
\begin{equation}\label{eq-solicitada-revisor}
\lambda(\mu_{X'})=\frac{1}{P'}\log|T'_{\bar{\omega}}(x')|=\frac{\log|T'_{\bar{\omega}}(x')|}{nP+|\bar\kappa_1|+|\bar\kappa_2|}.
\end{equation}
Noting that, $\log|T'_{\bar{\omega}}(x')|<0$ since $T_{\bar{\omega}}$ is a contraction, and substituting \eqref{iii.37} into equation~\eqref{eq-solicitada-revisor}, we have
\begin{equation}\label{iii.38}
 \begin{aligned}
  \lambda(\mu_{X'})>\frac{1}{nP} \left(C_{3}+n(1-\frac{2\log\nu}{3\log L_{1}})\log\alpha\right)
  =\lambda(\mu_X)(1-\frac{2\log\nu}{3\log L_{1}})+\frac{C_3}{nP}~,
 \end{aligned}
\end{equation}
where the last equality follows since $\lambda(\mu_X)=(\log \alpha)/P$. Since $C_3=C_2+C_1$ and the constants $C_1$ and $C_2$ in~\eqref{eq:C1} and~\eqref{eq:C2} do not depend on $n$ (they only depend on $\mathscr{F}$ and the interval $J$), then $C_3/nP \to 0$ as $n\to \infty$. Thus, similarly as above, since
$$
    \lambda(\mu_X)(1-\frac{\log\nu}{2\log L_{1}})<\lambda(\mu_X)(1-\frac{2\log\nu}{3\log L_{1}}),
$$
we find a sufficiently large $n$ such that~\eqref{iii.38} becomes
$$
  \lambda(\mu_{X'})>\lambda(\mu_X)(1-\frac{2\log\nu}{3\log
  L_{1}})+\frac{C_3}{nP}>\lambda(\mu_X)(1-\frac{\log\nu}{2\log L_{1}}).
$$
This proves Item~\ref{item1-lema3.1} of the lemma with $c=1-\frac{\log\nu}{2\log L_{1}}\in (0,1)$.

Now we will prove that $X'$ is a $(\aleph,\gamma)$-good approximation of  $X$. For this, $X'$ must satisfy conditions \ref{item-boa1}--\ref{item-boa3} in Definition~\ref{defi-boa-aprox}. Indeed, let $(\omega',x')$ and $(\omega,x)$ be periodic points of $F$ that generate $X'$ and $X$ respectively. 
Note that $\omega'=\bar{\omega}^\infty$ and $\omega=\hat{\omega}^\infty$. Since $\bar{\omega}=\hat{\omega}\,\overset{n}{\ldots}\, \hat{\omega}\bar\kappa_1\bar\kappa_2$ and $|\hat{\omega}|=P$, the first $Pn$ letters of the sequences  $\omega'$ and $\omega$ coincide. Let  $M$ be the smallest positive integer such that $2^{-MP}<\gamma$. Since $n$ can be arbitrarily large, we can also demand that $2M+1<n$. Then $M<n-M-1<n$, and therefore
\begin{equation}\label{iii.42}
 \tilde{X'}=\{F^{j}(\omega',x'): \, MP\leq j<(n-M-1)P \},
\end{equation}
is contained in the set of the first $n$ iterations by $F^P$ of the point $(\omega',x')$. Now let us define the projection $\rho:\tilde{X'}\to X$ as follows:
\begin{equation}\label{iii.43}
 \rho(F^{j}(\omega',x'))=F^{j}(\omega,x).
\end{equation}
Then, for $i=0,\dots,P-1$, $j=MP,\dots,(n-M-1)P-1$, we have
\begin{align*}
  d\left(F^i\left(F^j(\omega',x')\right), F^i\left(\rho\circ  F^j(\omega',x')\right)\right) 
  &= d\left(F^{i+j}(\omega',x'),F^{i+j}(\omega,x)\right) \\
 &= \max\{\bar{d}(\sigma^{i+j}(\omega'),\sigma^{i+j}(\omega)),
  |f^{i+j}_{\omega'}(x')-f^{i+j}_{\omega}(x)| \}.
\end{align*}
Since the first $nP$ letters of $\omega'$ and $\omega$ are the same, and since $MP\leq i+j<nP$, we have $\bar{d}(\sigma^{i+j}(\omega'),\sigma^{i+j}(\omega))\leq 2^{-MP}<\gamma$. For the same reason, writing $j=qP+\tau$, we have
$$
f^{i+j}_{\omega'}=f^{i+j}_{\omega}=f^i_{\sigma^j(\omega)}\circ
f^j_\omega = f^i_{\sigma^\tau(\omega)}\circ f^\tau_{\omega}\circ
f^q=f^{i+\tau}_\omega \circ f^q.
$$
Then, since $x$ and $x'$ are in $J$ and $f^q(J)\subset J$, by virtue of \eqref{iii.10}, it holds that
$
 |f^{i+j}_{\omega'}(x')-f^{i+j}_{\omega}(x)| \leq
 |f^{i+\tau}_\omega(J)|\leq L^{2P}|J|<\gamma$.
We conclude
$$
 d\left(F^i(y), F^i(\rho(y))\right)<\gamma \quad
 \text{for all $y\in \tilde{X'}$ and  $i=0,\dots,P-1$},
$$
proving Item~\ref{item-boa1} in Definition~\ref{defi-boa-aprox}.

Now let us show Item~\ref{item-boa2} in Definition~\ref{defi-boa-aprox}. For this, note that
\begin{equation}\label{iii.44}
\begin{aligned}
1-\frac{\#\tilde{X'}}{\#X'}&=\frac{|\bar{\kappa}_1|+|\bar\kappa_2|+(2M+1)P}{nP+|\bar{\kappa}_1|+|\bar\kappa_2|}
\leq \frac{|\bar{\kappa}_1|+K+(2M+1)P}{nP}. 
\end{aligned}
\end{equation}
On the other hand, denoting by $H=\max\{|\bar\beta_1|,\dots,|\bar\beta_m|\}$ where $\bar\beta_i$ is the 
finite word that generates $h_i$, 
by~\eqref{iii.15}, \eqref{eq:Pablo}, and $\bar\kappa_1=\bar\beta_{i_r}\dots\bar\beta_{i_1}$, we have
\begin{align*}
  |\bar{\kappa}_1| &\leq r H \leq \frac{-n\log \alpha_+ + \log \delta |J|^{-1}}{\log
  L_1} H \\
  &\leq \left(-n\log \alpha - K\log L\right) \frac{H}{\log L_1} =
  - \left (n {P  \, \lambda(\mu_X)} + {K\log L }\right) \frac{H}{\log
  L_1}.
\end{align*}
Substituting this expression into \eqref{iii.44} we have
$$
 1-\frac{\#\tilde{X'}}{\#X'} \leq |\lambda(\mu_X)| \frac{H}{\log
 L_1} + \frac{C_4}{nP}~,
$$
where  $C_4=  K + (2M+1)P-K H \frac{\log L}{\log L_1}$. Since $C_4$ does not depend on $n$, 
by again taking $n$ sufficiently large, it follows that
$$
\frac{C_{4}}{nP}+|\lambda(\mu_{X})|\frac{H}{\log L_1}\leq
|\lambda(\mu_{X})|\frac{2H}{\log L_1}.
$$
Therefore, from \eqref{iii.44}, we have
$$
\aleph=1-d\cdot |\lambda(\mu_X)|\leq\frac{\#\tilde{X'}}{\#X'} \quad
\text{where \( d=2H/\log L_1>0 \)}.
$$
This shows Item~\ref{item-boa2} in Definition~\ref{defi-boa-aprox}.

Finally, we only need to show Item~\ref{item-boa3} in Definition~\ref{defi-boa-aprox}. For this, given $\tau=0,1,\dots,P-1$, we have
\begin{equation*}\label{iii.45}
\begin{aligned}
\#\rho^{-1}(F^{\tau}(\omega,x))=&\#\{F^{j}(\omega',x')\in\tilde{X}': \, \rho(F^{j}(\omega',x'))=F^{\tau}(\omega,x)\}\\
=&\#\{\, j \in \mathbb{N} : MP \le j < (n - M - 1)P \text{ and } j= qP+\tau\, \text{for some}\, q\in\mathbb{N}\} \\
=& \frac{(n-M-1)P-MP}{P}=n-2M-1.
\end{aligned}
\end{equation*}
This proves the last item and therefore $X'$ is a $(\gamma,\aleph)$-good approximation of $X$.
\end{proof}

\subsection{Proof of Theorem~\ref{teo-B}} 


To prove this result, we use Theorem~\ref{teoBonati}. First, we establish that $F$ admits a sequence of periodic orbits $\{X_n\}$ with increasing periods $P_n$.

Let $\{\gamma_n\}_{n \in \mathbb{N}}$ and $\{\eta_n\}_{n \in \mathbb{N}}$ be positive sequences satisfying:
\begin{equation}\label{iii.46}
\sum_{n=1}^\infty \gamma_n < \infty \quad \text{and} \quad \eta_n \to 0 \quad \text{as} \quad n \to \infty.
\end{equation}
Let $X_1$ be the periodic orbit of $F$ generated by $f$, with period $P_1$. By recursively applying Lemma~\ref{lema3.1tesis}, for each $\gamma_n > 0$ and $\eta_n > 0$, there exists a periodic orbit $X_{n+1}$ of $F$ with period $P_{n+1} > P_n$ such that
\begin{enumerate}[label=\alph*)]
    \item $c \cdot \lambda(\mu_{X_n}) < \lambda(\mu_{X_{n+1}}) < 0$;
    \item $X_{n+1}$ is a $(\gamma_n, \aleph_n)$-good approximation of $X_n$, with $\aleph_n = 1 - d \cdot |\lambda(\mu_{X_n})|$;
    \item\label{item(c)-teo-principal} $X_{n+1}$ is $\eta_n$-dense in $\Sigma_k \times \mathbb{S}^1$.
\end{enumerate}
Here, $0 < c < 1$ and $d > 0$ are constants depending only on $\mathscr{F}$. In particular,
\[
c^n \cdot \lambda(\mu_{X_1}) < \lambda(\mu_{X_{n+1}}) < 0,
\]
which implies that $\lambda(\mu_{X_{n+1}}) \to 0$ as $n \to \infty$. Therefore, $\aleph_n > 0$ for sufficiently large $n$. Without loss of generality, we assume this holds for all $n \geq 1$. 

Furthermore, since
\[
\sum_{n=1}^\infty |\lambda(\mu_{X_n})| \leq \sum_{n=0}^\infty c^n |\lambda(\mu_{X_1})| < \infty,
\]
the Comparison Test for Series implies
\[
-\infty < \sum_{n=1}^\infty \log \left(1 - d \cdot |\lambda(\mu_{X_n})|\right) = \sum_{n=1}^\infty \log \aleph_n < 0,
\]
and consequently,
$\prod_{n=1}^\infty \aleph_n \in (0, 1]$.

By Theorem~\ref{teoBonati}, the sequence $\{\mu_{X_n}\}$ of atomic measures uniformly distributed over the orbits $X_n$ has a limit. The limit measure $\mu$ is ergodic, and its topological support is
$
\mathrm{supp}~\mu = \bigcap_{j=1}^\infty \left(\overline{\bigcup_{l=j}^\infty X_l}\right)$.
From \ref{item(c)-teo-principal} and the fact that $\eta_n \to 0$ as $n \to \infty$, we deduce that 
$
\overline{\bigcup_{l=j}^\infty X_l} = \Sigma_k \times \mathbb{S}^1$,
which implies $\mathrm{supp}~\mu = \Sigma_k \times \mathbb{S}^1$.

On the other hand, since the sequence of Lyapunov exponents $\{\lambda(\mu_{X_n})\}$ converges to zero, Proposition~\ref{pro3.3-tesis} implies that $\lambda(\mu) = 0$. This completes the proof of the theorem.

\section{Proof of Theorem~\ref{thmAA}} \label{s-thmA}
  For completeness, we recall that a family $\mathscr{F} \in \mathrm{IFS}^r_k(\mathbb{S}^1)$ is said to be \emph{transitive} if the associated skew product $F = \sigma \ltimes \mathscr{F}$ is transitive. Equivalent definitions can be found in~\cite[Def.~4.1, Thm.~4.2, 4.5, and 6.1]{BC23}. Moreover, $\mathscr{F}$ is called \emph{robustly minimal (transitive)} if any sufficiently small $C^1$ perturbation of $\mathscr{F}$ remains minimal (transitive).
  
  It is straightforward to verify that properties~\ref{item2-teo-B} and~\ref{item3-teo-B} of Theorem~\ref{teo-B} are open in the $C^1$ topology. This implies that any IFS sufficiently close to $\mathscr{F}$ in the $C^1$ topology also satisfies these properties. However, property~\ref{item1-teo-B} is not generally open. Nevertheless, it was established in~\cite{BFMS16} that if $\mathscr{F}$ is both expanding and minimal, then minimality is $C^1$-robust. Since minimality of $\mathscr{F}$ implies transitivity of the skew product $F=\sigma\ltimes \mathscr{F}$ (cf.~\cite[Remark~4.4 and Theorem~6.1]{BC23}), we obtain the following observation:

\begin{obs} \label{rem-1} 
The set of skew products $F = \sigma \ltimes \mathscr{F} \in S^r_k(\mathbb{S}^1)$ such that $\mathscr{F}$ satisfies properties \ref{item1-teo-B}--\ref{item3-teo-B} and 
\begin{itemize}
    \item[\mylabel{item4-obs}{(iv)}] $\mathscr{F}$ is forward expanding,
\end{itemize}
is an open subset of $\mathrm{TS}^r_k(\mathbb{S}^1)$ for every $r \geq 1$.
\end{obs}

\begin{defi} \label{def:globalized-blending-regions}
An open set $B$ is called a \emph{globalized blending region for $\mathscr{F}$} if there exist maps 
$h_1, \dots, h_m, T_1, \dots, T_s, S_1, \dots, S_t \in \langle \mathscr{F}\rangle^+$ 
and an open set
$D \subset \mathbb{S}^1$ such that $\overline{B} \subset D$ and
\begin{enumerate}[label=(\arabic*)]
    \item\label{item1-defi-region-misturadora} $\overline{B} \subset h_1(B) \cup \dots \cup h_m(B)$, 
    \item\label{item2-defi-region-misturadora} $h_i: \overline{D} \to D$ is a contracting map for $i=1, \dots, m$,
    \item\label{item3-defi-region-misturadora} $\mathbb{S}^1 = T_1(B) \cup \dots \cup T_s(B) = S_1^{-1}(B) \cup \dots \cup S_t^{-1}(B)$.
\end{enumerate}
\end{defi}

As shown in~\cite[Sec.~6]{BFMS16}, a skew product $F = \sigma \ltimes \mathscr{F}$ where both $\mathscr{F}$ and $\mathscr{F}^{-1}$ possess globalized blending regions satisfies conditions~\ref{item1-teo-B}--\ref{item4-obs}.  In particular, Theorem~\ref{teo-B} applies. 
This result parallels~\cite[Prop.~8.3]{diaz2017nonhyperbolic}, where robust transitivity follows from the existence of blenders (analogous to conditions~\ref{item1-defi-region-misturadora} and \ref{item2-defi-region-misturadora} of Definition~\ref{def:globalized-blending-regions}) combined with the requirement that every point has forward and backward iterates inside the blender's domain (analogous to condition~\ref{item3-defi-region-misturadora}).

The key idea behind the proof of Theorem~\ref{thmAA} is to demonstrate that, up to conjugation, any transitive locally constant skew product can be approximated by one satisfying conditions~\ref{item1-teo-B}--\ref{item4-obs}. These properties essentially arise from constructing a pair of globalized blending regions for $\mathscr{F}$ and $\mathscr{F}^{-1}$, which can be found within the proof of the next results borrowed from~\cite{BC23}.

The following result was proved by Barrientos and Cisneros in~\cite[Cor.~5.5]{BC23}.  

\begin{teo}\label{teo-aproxi-elem-D}
Let $\mathscr{F}$ be a transitive family in $\IFS^r_k(\mathbb{S}^1)$, where $r \geq 1$ and $k \geq 2$. Then, there exists $\mathscr{G} \in \IFS^r_k(\mathbb{S}^1)$, $C^r$-close to $\mathscr{F}$, such that $\mathscr{G}$ is both robustly forward and backward minimal, the elements of $\mathscr{G}$ are Morse-Smale diffeomorphisms, and there exists a map $g_0 \in \langle \mathscr{G}\rangle^+$ that is $C^\infty$-conjugate to an irrational rotation.  
\end{teo}

The next result, stated in~\cite[Thm.~3.25]{BC23}, improves~\cite[Thm.~B]{BFS14} and generalizes earlier results on robustly minimal IFSs on the circle described in~\cite{GI99, GI00}.

\begin{teo}\label{teo-B-BFS14}
Let $f_1, f_2$ be orientation-preserving circle diffeomorphisms. Assume that $f_1$ is an irrational rotation and $f_2$ is not a rotation. Then, $\mathscr{F} = \{f_1, f_2\}$ is both forward and backward expanding and robustly forward and backward minimal.
\end{teo}

Now we are ready to give the proof of our main result. 

\begin{proof}[Proof of Theorem~\ref{thmAA}] 
Denote by $\mathrm{T}^r_k(\mathbb{S}^1)$ the set of transitive elements of $\mathrm{IFS}^r_k(\mathbb{S}^1)$. Since every minimal IFS is also transitive (cf.~\cite[Rem.~4.4]{BC23}), by Theorem~\ref{teo-aproxi-elem-D}, there is a dense set $\mathcal{D}$ of $\mathrm{T}^r_k(\mathbb{S}^1)$ such that every $\mathscr{G}\in \mathcal{D}$ is robustly transitive, its elements are Morse-Smale diffeomorphisms, and there is a map $g_0\in \langle \mathscr{G}\rangle^+$ that is $C^\infty$ conjugate to an irrational rotation.    
Let $\phi$ be the $C^\infty$ conjugacy map between $g_0$ and the rigid irrational rotation. Set $\hat{\mathscr{G}}=\phi\circ \mathscr{G}\circ\phi^{-1}$. Note that since the property of being a Morse-Smale diffeomorphism is invariant under smooth conjugation, the elements of $\hat{\mathscr{G}}$ are also Morse-Smale diffeomorphisms. Moreover, since $\langle\hat{\mathscr{G}}\rangle^+$ also has an irrational rotation, by Theorem~\ref{teo-B-BFS14}, we conclude that $\hat{\mathscr{G}}$ is both backward and forward expanding and robustly forward and backward minimal. Since the expanding
property is also a robust property (cf.~\cite[Cor.~3.1]{BFMS16} or \cite[Rem.~3.10]{BC23}) and $\mathscr{G}$ is robustly transitive,  we can consider a neighborhood $\mathcal{U}_{\mathscr{G}} \subset \mathrm{T}^r_k(\mathbb{S}^1)$ of $\mathscr{G}$ in $\IFS^r_k(\mathbb{S}^1)$ such that for every $\mathscr{H}\in \mathcal{U}_{\mathscr{G}}$, the family $\hat{\mathscr{H}}=\phi\circ \mathscr{H}\circ\phi^{-1}$ is both backward and forward expanding and robustly forward and backward minimal. Thus, by Theorem~\ref{teo-B}, the skew product $\hat{H}=\sigma\ltimes\hat{\mathscr{H}}$ has an $\hat H$-invariant ergodic probability measure $\hat\mu$ with full support, limit of periodic measures and with Lyapunov exponent $\lambda(\hat\mu)=0$. Consider the map $\Phi=\mathrm{id}_{\Sigma_{k}}\times\phi$ and notice that $\hat{H}=\Phi\circ H\circ \Phi^{-1}$, where $H=\sigma\ltimes\mathscr{H}$. 

\begin{claim} \label{claim-final} The measure $\mu\eqdef\hat\mu \circ \Phi$ is an $ H$-invariant ergodic probability with full support, limit of periodic measures, and with Lyapunov exponent $\lambda(\mu)=0$.
\end{claim}
\begin{proof} 
We first observe that $\hat{H}^{-1} = \Phi \circ H^{-1} \circ \Phi^{-1}$. Hence, for any measurable set $A \subset \Sigma_{k} \times \mathbb{S}^{1}$, since $\hat{\mu}$ is a $\hat{H}$-invariant measure, i.e., $(\hat{H})_{*}\hat{\mu} = \hat{\mu}$, we have 
\begin{equation*}
    \mu(A) = \hat{\mu}(\Phi(A)) = (\hat{H})_{*}\hat{\mu}(\Phi(A)) = \hat{\mu}(\Phi \circ H^{-1} \circ \Phi^{-1} \circ \Phi(A)) = \hat{\mu}(\Phi \circ H^{-1}(A)) = \mu(H^{-1}(A)).    
\end{equation*}
Thus, $\mu$ is $H$-invariant. Moreover, if $A$ satisfies $H^{-1}(A) = A$, then 
\begin{equation*}
    \hat{H}^{-1}(\Phi(A)) = \Phi \circ H^{-1} \circ \Phi^{-1}(\Phi(A)) = \Phi \circ H^{-1}(A) = \Phi(A)
\end{equation*}
and thus, since $\hat{\mu}$ is an ergodic measure for $\hat{H}$, we have that $\mu(A) = \hat{\mu}(\Phi(A)) \in \{0,1\}$. Hence, $\mu$ is also an ergodic measure for $H$. Also, since $\Phi$ is a homeomorphism and $\hat{\mu}$ has full support, then $\mu = \hat{\mu} \circ \Phi$ inherits full support. Furthermore, $\mu$ is also a limit of periodic measures, as $\hat{\mu}$ is.
Finally, we will calculate the Lyapunov exponent $\lambda(\mu)$. To do this, observe that $\hat{\mathscr{H}}=\phi\circ \mathscr{H} \circ \phi^{-1}$ and thus the fiber maps of $\hat{\mathscr{H}}^n$  are given by $\hat{h}^n_\omega=\phi \circ h^n_\omega \circ \phi^{-1}$  for any $\omega\in \Sigma_k$, where $h^n_\omega$ are the fiber maps of $H^n$. Hence,
\begin{equation*}
\begin{aligned}
     \lambda_{\hat{H}}(\omega,x) 
     & = \lim_{n\to\infty} \frac{1}{n} \log |(\hat h^n_\omega)'(x)| \\ 
     &=  \lim_{n\to\infty} \frac{1}{n} \log |\phi'(h^n_\omega\circ \phi^{-1}(x))|+ \lim_{n\to\infty} \frac{1}{n} \log |(h^n_\omega)'(\phi^{-1}(x))| + \lim_{n\to\infty} \frac{1}{n} \log |(\phi^{-1})'(x)| \\
     &= \lim_{n\to\infty} \frac{1}{n} \log |(h^n_\omega)'(\phi^{-1}(x))| = \lambda_{H}(\omega,\phi^{-1}(x))= \lambda_{H}(\Phi^{-1}(\omega,x))
\end{aligned}
\end{equation*}
for every $(\omega,x)\in \Sigma_k\times\mathbb{S}^1$.  Now, since $\hat\mu= \mu \circ \Phi^{-1}$, by the above equality,  we have
$$
\lambda(\hat{\mu})=\int \lambda_{\hat{H}}(\omega,x) \, d\hat\mu = \int \lambda_{\hat{H}}(\Phi^{-1}(\omega,x)) \, d\mu =\int \lambda_{{H}}(\omega,x) \, d\mu = \lambda(\mu).
$$
Since $\lambda(\hat\mu)=0$, we also get that $\lambda(\mu)=0$ as required. This completes the proof
\end{proof}

Now, we conclude the proof by taking $\tilde{\mathcal{R}}$ as the open set formed by the union of the neighborhoods $\mathcal{U}_{\mathscr{G}}$ for $\mathscr{G}\in \mathcal{D}$. Since $\mathcal{U}_{\mathscr{G}}\subset \mathrm{T}^r_{k}(\mathbb{S}^1)$ and $\mathcal{D}$ is dense in $\mathrm{T}^r_{k}(\mathbb{S}^1)$, we get that $\tilde{\mathcal{R}}$ is open and dense in $\mathrm{T}^r_{k}(\mathbb{S}^1)$. Then, by the definition of the topology in $\mathrm{S}^r_k(\mathbb{S}^1)$ and in view that $\mathscr{F}\in \mathrm{T}^r_{k}(\mathbb{S}^1)$  if and only if $F=\sigma\ltimes \mathscr{F}\in \mathrm{T}^r_{k}(\mathbb{S}^1)$, we get the set $\mathcal{R}=\{F=\sigma\ltimes \mathscr{F}\in \mathrm{S}^r_{k}(\mathbb{S}^1): \mathscr{F} \in \tilde{\mathcal{R}}\}$ is an open and dense subset of $\mathrm{TS}^r_k(\mathbb{S}^1)$. Moreover, for every $H=\sigma\ltimes \mathscr{H}\in \mathcal{R}$, the family $\mathscr{H}$ belongs to $\mathcal{U}_{\mathscr{G}}$ for some $\mathscr{G}\in \mathcal{D}$ and thus by Claim~\ref{claim-final}, there is an ergodic $F$-invariant probability measure $\mu$ with full support whose Lyapunov exponent $\lambda(\mu)$ equal to zero. Furthermore, the probability measure $\mu$ is a limit in the weak$^*$ topology of sequences of invariant measures supported on periodic points. This concludes the proof.
\end{proof}

\section*{Acknowledgement}
The first author, P.~G.~Barrientos, was supported by grant PID2020-113052GB-I00 funded by MCIN, PQ 305352/2020-2 (CNPq), and JCNE E-26/201.305/2022 (FAPERJ). The proof of Theorem~\ref{teo-B} forms part of the results obtained in the master's thesis of the author, J.~A.~Cisneros. We express our deepest gratitude to Prof.~Lorenzo J.~D\'iaz for his insightful comments on the results obtained. His contributions as a member of the master's and doctoral committees of the first author have been invaluable and greatly appreciated.

\newcommand{\etalchar}[1]{$^{#1}$}

\end{document}